\documentclass[12pt]{article}

\usepackage{amssymb,amsmath,amsthm}

\newcommand{\ovl}{\overline}
\newcommand{\n}{\noindent}
\newcommand{\vp}{\varepsilon}
\newcommand{\cl}[1]{{\mathcal{#1}}}
\newcommand{\bb}[1]{{\mathbb{#1}}}

\theoremstyle{plain}
\newtheorem{thm}{Theorem}[section]

\newtheorem{lem}[thm]{Lemma}
\newtheorem{cor}[thm]{Corollary}

\theoremstyle{definition}
\newtheorem{defn}[thm]{Definition}

\theoremstyle{remark}
\newtheorem{rem}[thm]{Remark}

\begin{document}

\title
{THE LAPLACIAN MASA IN \\A FREE GROUP FACTOR}
\author{Allan M.\ Sinclair\\
Department of Mathematics\\
University of Edinburgh\\
Edinburgh EH9 3JZ\\
SCOTLAND\\
{\tt e-mail:\ allan@maths.ed.ac.uk}
\and
Roger R. Smith\footnote
{Partially supported by a grant from the National Science  
Foundation. }\\
Department of Mathematics\\
Texas A\&M University\\
College Station, TX \ 77843\\ 
U.S.A.\\
{\tt e-mail:\ rsmith@math.tamu.edu}}

\date{}

\maketitle{}

\bigskip

\begin{abstract}{
The Laplacian (or radial) masa in a free group factor is generated by the sum
of the generators and their inverses. We show that such a masa $\mathcal{B}$
is strongly singular and has Popa invariant $\delta(\cl B) = 1$. This is
achieved by proving that the conditional expectation $\bb E_{\cl B}$ onto $\cl
B$ is an asymptotic homomorphism. We also obtain similar results for the free
product of discrete groups, each of which contains an element of infinite
order. }\end{abstract}
\hfill

 \newpage 
\setcounter{section}{0}

\section{Introduction}\label{sec1}

\indent

A maximal abelian self-adjoint subalgebra $\cl A$ of a type ${\rm II}_1$ factor
$\cl M$ is singular, \cite{Dix}, if every unitary $u\in \cl M$ which normalizes
$\cl A$ must lie in $\cl A$. Popa, \cite{Po1983}, has shown that every type
${\rm II}_1$ factor contains such a masa, and there are explicit examples in
various group von~Neumann factors, \cite{BR,Dix,Po1983,Ra,SS2001}. Several
invariants have been introduced for singular masas in type ${\rm II}_1$
factors, of which Popa's delta invariant, \cite{Po1983}, strong singularity
(see (2.3) and \cite{SS2001}) and the existence of an asymptotic homomorphism,
\cite{SS2001}, are relevant for this paper. We prove that the Laplacian, or
radial, masa $\cl B$ in a free group factor is strongly singular, has
$\delta(\cl B) = 1$, and its associated conditional expectation $\bb E_{\cl
B}$ is an asymptotic homomorphism, \cite[Section~5]{SS2001}.

Let $g_1,\ldots, g_k$ be the generators of the free group $\bb F_k$, $k\ge 2$,
and let
$$w_1 = \sum^k_{i=1} (g_i+g^{-1}_i).$$
Then $\cl B$ is defined to be the abelian subalgebra of $VN(\bb F_k)$
generated by the self-adjoint element $w_1$. Pytlik, \cite{Py}, showed that
$\cl B$ is a masa, and subsequently R\u{a}dulescu, \cite{Ra}, proved that it is
singular by establishing that the Puk\'{a}nszky invariant, \cite{Pu}, is
$\{1,\infty\}$. Popa, \cite{Po1985}, had already demonstrated that a masa with
this particular Puk\'{a}nszky invariant is singular.

Our aim in this paper is to show that $\bb E_{\cl B}$ is an asymptotic
homomorphism, Theorem~\ref{thm5.1}, since strong singularity and $\delta(\cl
B) = 1$ are then implied by \cite[Theorem~6.7]{SS2001}. We achieve this by
proving a sufficient condition, (Theorem~\ref{thm3.1}), phrased in terms of
the convergence of certain infinite series. Specifically, we need that the
difference between $\bb E_{\cl B}(xw_ny)$ and $\bb E_{\cl B}(x) \bb E_{\cl
B}(y) w_n$, where $w_n$ is the sum of words of length $n$ and $x$ and $y$ are
fixed words in $\bb F_k$, should tend to 0 rapidly as $n\to \infty$.
We also show, in Theorem \ref{thm5.4}, that our results extend to certain free
products of groups (see \cite{BR} for similar work).

The second section contains background material and definitions of the terms
mentioned above. The heart of the paper is the fourth section in which we
obtain various enumeration results for words in $\bb F_k$. These are used to
show that the Laplacian masa satisfies the hypothesis of Theorem~\ref{thm3.1},
from which strong singularity follows.\newpage

\section{Preliminaries}\label{sec2}

\indent

In this section we present definitions and notation which will be needed
subsequently. Specific notation about the number of words starting and ending
in a certain way is contained in the fourth section, where it is used.

The operator norm on a type ${\rm II}_1$ factor $\cl M$ with normalized trace
tr is denoted by $\|\cdot\|$, and $\|\cdot\|_2$ is the norm $\|x\|_2 =
(\text{tr}(x^*x))^{1/2}$. If $\cl N$ is a von~Neumann subalgebra of $\cl M$
then $\bb E_{\cl N}$ denotes the trace preserving conditional expectation of
$\cl M$ onto $\cl N$. If $\phi\colon \ \cl M_1\to \cl M_2$ is a linear map
between type ${\rm II}_1$ factors, then there are several norms on $\phi$,
depending on the norms given to $\cl M_1$ and $\cl M_2$. We write $\|\phi\|$
and $\|\phi\|_2$ when both factors have respectively the operator norm and the
$\|\cdot\|_2$-norm. When $\cl M_1$ has the operator norm and $\cl M_2$ has the
$\|\cdot\|_2$-norm, we denote the resulting norm of $\phi$ by
$\|\phi\|_{\infty,2}$. This last norm was introduced and studied in
\cite{SS2001}. We note that conditional expectations are contractions for all
three norms.

A maximal abelian self-adjoint subalgebra (masa) in $\cl M$ is singular if any
unitary $u\in \cl M$ which normalizes $\cl A$, $(u\cl A u^* = \cl A)$, must
lie in $\cl A$, \cite{Dix}. A masa $\cl A$ has Popa invariant $\delta(\cl A) =
1$ if and only if the following is satisfied: \ for each nilpotent $(v^2=0)$
partial isometry $v\in\cl M$ with $v\cl Av^*$ and $v^*\cl Av$ contained in
$\cl A$, we have
\begin{equation}\label{eq2.1}
\sup\{\|(I-\bb E_{\cl A})(x)\|_2\colon \ x\in v\cl Av^*, \|x\|\le 1\}\ge
\|vv^*\|_2,
\end{equation}
\cite{Po1983,SS2001}. In \cite{SS2001} we introduced the concept of strong
singularity for masas $\cl A$. The defining property is
\begin{equation}\label{eq2.2}
\|u-\bb E_{\cl A}(u)\|_2 \le \|\bb E_{u\cl Au^*}- \bb E_{\cl A}\|_{\infty,2}
\end{equation}
for all unitaries $u\in \cl M$. Singularity is clearly implied by
(\ref{eq2.2}), and we note that a reverse inequality
\begin{equation}\label{eq2.3}
\|\bb E_{u\cl Au^*}  -\bb E_{\cl A}\|_{\infty,2} \le 4\|u-\bb E_{\cl A}(u)\|_2.
\end{equation}
holds for all masas, \cite{SS2001}. In the same paper, we introduced the notion
of an asymptotic homomorphism. This is the conditional expectation $\bb
E_{\cl A}$ onto an abelian subalgebra $\cl A$ for which a unitary $u\in \cl A$
may be found to satisfy
\begin{equation}\label{eq2.4}
\lim_{|k|\to \infty} \|\bb E_{\cl A}(xu^ky) - \bb E_{\cl A}(x) \bb E_{\cl
A}(y) u^k\|_2 = 0
\end{equation}
for all $x,y\in \cl M$. The masa associated to a generator in the von~Neumann
algebra of a free group  is one particular example of this phenomenon,
\cite{SS2001}.

The above are all rigidity conditions on the abelian von~Neumann subalgebra
$\cl A$ of $\cl M$. The main point for asymptotic homomorphisms is that their
existence implies both singularity and strong singularity, and that
$\delta(\cl A) = 1$, \cite[Theorem~6.7]{SS2001}.

Let $\bb F_k$ denote the free group on $k$ generators, $2\le k < \infty$. We
denote the generators of $\bb F_k$ by $\{g_i\}^k_{i=1}$, and we let
$$S_k = \{g_1,\ldots, g_k,g^{-1}_1,\ldots, g^{-1}_k\}$$
with cardinality $|S_k|= 2k$. For $n\ge 0$, we define $w_n$ to be the sum of
all reduced words in $\bb F_k$ of length $n$. Then, for $n\ge 2$,
\begin{align}
w_0 &= e,\nonumber\\
w_1 &= \sum^k_{i=1} g_i + \sum^k_{i=1} g^{-1}_i,\nonumber\\
w^2_1 &= w_2 + 2kw_0,\nonumber\\
\label{eq2.5}
w_1w_n &= w_nw_1 = w_{n+1} + (2k-1)w_{n-1}
\end{align}
are examples of relations among these sums. They also show that the unital
algebra generated by $w_1$ contains each $w_n$, and that this algebra is the
span of these elements. The words in $w_n$ are pairwise orthogonal with
respect to the trace, and so a simple counting argument gives the well known
formula
\begin{equation}\label{eq2.6}
\|w_n\|^2_2 = 2k(2k-1)^{n-1},\qquad n\ge 1,
\end{equation}
since this is the number of words of length $n$.

Let $\cl B$ denote the abelian von~Neumann algebra generated by $w_1$. Then
$\cl B$ is a masa, \cite{Py}, and our previous remarks show that
$$\{w_n/\|w_n\|_2 \colon \ n\ge 0\}$$
is an orthonormal basis for $L^2(\cl B, \text{ tr})$. Moreover, for each $x\in
\cl M$,
\begin{equation}\label{eq2.7}
\bb E_{\cl B}(x) = \sum^\infty_{n=0} \text{ tr}(xw_n)w_n/\|w_n\|^2_2.
\end{equation}\newpage

\section{Asymptotic homomorphisms}\label{sec3}

\setcounter{equation}{0}

\indent

In this section we present a criterion for determining when a conditional
expectation is an asymptotic homomorphism, phrased in terms of an orthonormal
basis for $L^2(\cl B, \text{ tr})$.

\begin{thm}\label{thm3.1}
Let $\cl A$ be an abelian von~Neumann subalgebra of a type ${\rm II}_1$ factor
$\cl M$, and suppose that there is a $*$-isomorphism $\pi\colon \ \cl A\to
L^\infty[0,1]$ which induces an isometry from $L^2(\cl A, \text{ tr})$ onto
$L^2[0,1]$. Let
$$\{v_n\colon \ n\ge 0,\quad v_n\in \cl A\}$$
be an orthonormal basis for $L^2(\cl A, \text{ tr})$, and let $Y\subseteq \cl
M$ be a set whose linear span is norm dense is $L^2(\cl M, \text{ tr})$. If 
\begin{equation}\label{eq3.1}
\sum^\infty_{n=0} \|\bb E_{\cl A}(xv_ny) - \bb E_{\cl A}(x) \bb E_{\cl A}(y)
v_n\|^2_2 < \infty
\end{equation}
for all $x,y\in Y$, then $\bb E_{\cl A}$ is an asymptotic homomorphism, $\cl
A$ is a strongly singular masa, and $\delta(\cl A) = 1$.
\end{thm}

\begin{proof}
Suppose that there is a unitary $u\in \cl A$ for which
\begin{equation}\label{eq3.2}
\lim_{|n|\to\infty}\|\bb E_{\cl A}(xu^ny) - \bb E_{\cl A}(x) \bb E_{\cl A}(y)
u^n\|_2 = 0
\end{equation}
when $x$ and $y$ are arbitrary elements of $Y$. Then (\ref{eq3.2}) will hold
for $x,y\in \text{span } Y$, and the $\|\cdot\|_2$-norm continuity of $\bb
E_{\cl A}$ will extend its validity to all $x,y\in \cl M$. Thus we may
restrict attention to $x,y\in Y$.

Let $w$ be the unitary $e^{2\pi it}$, $0 \le t\le 1$, which generates
$L^\infty[0,1]$, and let $u = \pi^{-1}(w)$ be the corresponding unitary
generator of $\cl A$. Since $\pi$ is an isometry for the Hilbert space norms,
we have
\begin{equation}\label{eq3.3}
\text{tr}(v^*_ju^n) = \langle u^n,v_j\rangle= \int^1_0 \,e^{2\pi int}
\ovl{\pi(v_j)(t)} \ dt.
\end{equation}
From this it follows that
\begin{equation}\label{eq3.4}
\lim_{|n|\to\infty} \langle u^n,v_j\rangle = 0,\qquad j\ge 0,
\end{equation}
because the right hand side of (\ref{eq3.3}) is the conjugate of the 
$n^{\rm
th}$ Fourier coefficient of $\pi(v_j) \in L^2[0,1]$.

Fix $x,y\in Y$, let $\vp>0$, and choose $k$, by (\ref{eq3.1}), so that
\begin{equation}\label{eq3.5}
\sum_{j>k} \|\bb E_{\cl A}(xv_jy) - \bb E_{\cl A}(x) \bb E_{\cl A}(y)v_j\|^2_2
< \vp^2.
\end{equation}
From (\ref{eq3.4}), we may now choose $n_0$ so that
\begin{equation}\label{eq3.6}
\langle u^n,v_j\rangle \leq \vp((k+1)\|x\| \ \|y\|)^{-1},\qquad 0\le j \le k,
\quad |n|\ge n_0.
\end{equation}
Since $\{v_j\}^\infty_{j=1}$ is an orthonormal basis, we may write $u^n =
\sum\limits^\infty_{j=0} \langle u^n, v_j\rangle v_j$, which gives 
\begin{equation}\label{eq3.7}
\|\bb E_{\cl A}(xu^ny) - \bb E_{\cl A}(x) \bb E_{\cl A}(y)u^n\|_2 =
\left\|\sum^\infty_{j=0} \langle u^n,v_j\rangle (\bb E_{\cl A}(xv_jy) - \bb
E_{\cl A}(x) \bb E_{\cl A}(y)v_j)\right\|_2.
\end{equation}
We split this sum at $k$ and estimate each part separately. By (\ref{eq3.6}),
\begin{align}
&\left\|\sum^k_{j=0} \langle u^n,v_j\rangle (\bb E_{\cl A}(xv_jy) - \bb E_{\cl
A}(x) \bb E_{\cl A}(y)v_j\right\|_2\nonumber\\
&\quad \le \sum^k_{j=0} \vp((k+1)\|x\| \ \|y\|)^{-1} \|\bb E_{\cl A}(xv_jy)  -
\bb E_{\cl A}(x) \bb E_{\cl A}(y)v_j\|_2\nonumber\\
&\quad \le \sum^k_{j=0} \vp((k+1)\|x\| \ \|y\|)^{-1} 2\|x\| \ \|y\|\nonumber\\
\label{eq3.8}
&\quad = 2\vp,
\end{align}
provided that $|n|\ge n_0$. Since
\begin{equation}\label{eq3.9}
1 = \|u^n\|^2_2 = \sum^\infty_{j=0} |\langle u^n,v_j\rangle|^2,
\end{equation}
we may use the Cauchy-Schwarz inequality to estimate
\begin{align}
&\left\|\sum^\infty_{j=k+1} \langle u^n,v_j\rangle (\bb E_{\cl A}(xv_jy) - \bb
E_{\cl A}(x) \bb E_{\cl A}(y)v_j\right\|_2\nonumber\\
&\quad \le \sum^\infty_{j=k+1} |\langle u^n,v_j\rangle|\|\bb E_{\cl A}(xv_jy) -
\bb E_{\cl A}(x) \bb E_{\cl A}(y) v_j\|_2\nonumber\\
&\quad \le \left(\sum^\infty_{j=k+1} \|\bb E_{\cl A}(xv_jy) - \bb E_{\cl A}(x)
\bb E_{\cl A}(y) v_j\|^2_2\right)^{1/2}\nonumber\\
\label{eq3.10}
&\quad < \vp,
\end{align}
by the choice of $k$ (see (\ref{eq3.5})). Then (\ref{eq3.8}) and
(\ref{eq3.10}) allow us to obtain the inequality
\begin{equation}\label{eq3.11}
\|\bb E_{\cl A}(xu^ny) - \bb E_{\cl A}(x) \bb E_{\cl A}(y)u^n\|_2 <
3\vp,\qquad |n|\ge n_0,
\end{equation}
from (\ref{eq3.7}). This proves (\ref{eq3.2}), and the proof is complete.
\end{proof}\newpage

\section{Enumeration of words in $\pmb{\bb F_k}$}\label{sec4}

\setcounter{equation}{0}

\indent

This section contains the main technical results which we will use to show that
$\bb E_{\cl B}$ is an asymptotic homomorphism. We fix an integer $k\ge 2$, and
formulate our results for $\bb F_k$. Some of the quantities below depend on
$k$, but we have suppressed this for notational convenience.

\begin{defn}\label{defn4.1}
Let $\{g_1,\ldots, g_k\}$ be the generators of $\bb F_k$ and let
$$S_k = \{g_1,\ldots, g_k, g^{-1}_1,\ldots, g^{-1}_k\}.$$
For $x,y\in S_k$, we let $w_n(x,y)$ be the sum of all words of length $n$
beginning with $x$ and ending with $y$, and we denote by $\nu_n(x,y)$ the
number of such words. If $\sigma$ and $\tau$ are non--empty subsets of $S_k$,
we write $\nu_n(\sigma,\tau)$ for the number of words of length $n$ which
begin  with an element of $\sigma$ and end with an element of $\tau$.
$\hfill\square$
\end{defn}

We note that the formula
\begin{equation}\label{eq4.1}
\nu_n(\sigma,\tau) = \sum_{\stackrel{\scriptstyle x\in\sigma}{\scriptstyle
 y\in\tau}}\nu_n(x,y) 
\end{equation}
is an immediate consequence of these definitions, for any non--empty subsets
$\sigma$ and $\tau$ of $S_k$. The following lemma is obvious and we omit the
proof.

\begin{lem}\label{lem4.2}
Let $\sigma = S_k\backslash\{g_1,g_2,g^{-1}_1,g^{-1}_2\}$ and let $\tau = S_k
\backslash\{g_1,g^{-1}_1\}$. Then, for $n\ge 2$,
\begin{align}
\label{eq4.2}
w_{n+1}(g_1,g_2) &= g_1\left(w_n(g_1,g_2) + w_n(g_2,g_2) + w_n(g^{-1}_2,g_2) +
\sum_{x\in\sigma} w_n(x,g_2)\right);\\
\label{eq4.3}
w_{n+1}(g_1,g_1) &= g_1\left(w_n(g_1,g_1) + \sum_{x\in\tau}
w_n(x,g_1)\right);\\
\label{eq4.4}
w_{n+1}(g_1,g^{-1}_1) &=g_1\left(w_n(g_1,g^{-1}_1) + \sum_{x\in\tau}
w_n(x,g^{-1}_1)\right).
\end{align}
\end{lem}

For $n\ge 2$, we introduce three constants $\alpha_n$, $\beta_n$, and
$\gamma_n$ which are respectively $\nu_n(g_1,g_2)$, $\nu_n(g_1,g_1)$ and
$\nu_n(g_1,g^{-1}_1)$. For any pair $x,y\in S_k$, there is an automorphism of
$\bb F_k$  which takes $w_n(x,y)$ to one of $w_n(g_1,g_2)$, $w_n(g_1,g_1)$
and $w_n(g_1,g^{-1}_1)$, and $\nu_n(x,y)$ is $\alpha_n,\beta_n$ or $\gamma_n$.

\begin{lem}\label{lem4.3}
The following relations hold for $\alpha_n,\beta_n$ and $\gamma_n$:

\n (i)~~For $n=2$,
\begin{equation}\label{eq4.5}
\alpha_2 = 1,\quad \beta_2=1, \quad \gamma_2  =
0.\end{equation}

\n (ii)~~For $n \geq 2$,
\begin{align}\label{eq4.6}
\alpha_{n+1} &= (2k-3) \alpha_n + \beta_n +
\gamma_n, \\
\label{eq4.7}
\beta_{n+1} &= \beta_n + (2k-2) \alpha_n,\\
\label{eq4.8}
\gamma_{n+1} &= \gamma_n + (2k-2)\alpha_n.
\end{align}

\n (iii)~~There exists a constant $C_k$ such that, for $n\ge 2$,
\begin{equation}\label{eq4.9}
|\alpha_n - (2k-1)^{n-1}/2k|, ~~|\beta_n-(2k-1)^{n-1}/2k|,~~
|\gamma_n-(2k-1)^{n-1}/2k| \le C_k.
\end{equation}
\end{lem}

\begin{proof}

\n (i)~~Clear.

\n (ii)~~These are the result of counting the terms in (\ref{eq4.2}),
(\ref{eq4.3}) and (\ref{eq4.4}) and noting that the cardinalities of $\sigma$
and $\tau$ in Lemma~\ref{lem4.2} are respectively $2k-4$ and $2k-2$.

\n (iii)~~Equations (\ref{eq4.7}) and (\ref{eq4.8}) show that 
\begin{equation}\label{eq4.10}
\beta_{n+1} - \gamma_{n+1} = \beta_n - \gamma_n,
\end{equation}
and so these differences are independent of $n$. Since $\beta_2-\gamma_2=1$,
we conclude that
\begin{equation}\label{eq4.11}
\beta_n=1+\gamma_n,\qquad n\ge 2.
\end{equation}
Subtraction of (\ref{eq4.8}) from (\ref{eq4.6}) gives
\begin{equation}\label{eq4.12}
\alpha_{n+1} - \gamma_{n+1} = \beta_n-\alpha_n = 1 + (\gamma_n-\alpha_n),
\end{equation}
using (\ref{eq4.11}). Since $\alpha_2-\gamma_2=1$, a simple induction argument
based on (\ref{eq4.12}), shows that
\begin{equation}\label{eq4.13}
\alpha_n = \gamma_n + (1+(-1)^n)/2,\qquad n\ge 2.
\end{equation}
It follows from (\ref{eq4.11}) and (\ref{eq4.13}) that
\begin{equation}\label{eq4.14}
|\alpha_n-\gamma_n|\le 1,\quad |\alpha_n-\beta_n| \le 2,\qquad n\ge 2.
\end{equation}

There are $2k(2k-1)^{n-1}$ words in $w_n$, so counting these words according
to whether they lie in a sum of the form $w_n(a,b)$, $w_n(a,a)$ or
$w_n(a,a^{-1})$ leads to
\begin{equation}\label{eq4.15}
2k(2k-2) \alpha_n + 2k\beta_n + 2k\gamma_n = 2k(2k-1)^{n-1}
\end{equation}
which, after cancellation, is
\begin{equation}\label{eq4.16}
(2k-2)\alpha_n + \beta_n + \gamma_n = (2k-1)^{n-1}.
\end{equation}
Substituting from (\ref{eq4.14}) gives
\begin{equation}\label{eq4.17}
|2k\alpha_n-(2k-1)^{n-1}| \le 3.
\end{equation}
Now (\ref{eq4.9}) follows from (\ref{eq4.17}) and (\ref{eq4.14}) if we define
$C_k$ to be $2+3/2k$.
\end{proof}

\begin{rem}\label{rem4.4}
As is surely well known, equations (\ref{eq4.6})--(\ref{eq4.8}) give closed
form expressions for $\alpha_n$, $\beta_n$ and $\gamma_n$. One way to see this
is to rewrite these equations as $\xi_{n+1}=A\xi_n$, with solution 
$\xi_n=A^{n-2}\xi_2$, where
\[\xi_n=\left(\begin{array}{c}\alpha_n\\ \beta_n\\ \gamma_n \end{array}
\right),\ \ \ A=
\left(\begin{array}{ccc}2k-3&1&1\\2k-2&1&0\\2k-2&0&1\end{array}\right).\]
Explicit solutions can then be exhibited by observing that the eigenvalues of 
$A$ are $2k-1$, $1$ and $-1$, with respective eigenvectors
\[(1,1,1)^T,\ \ (0,1,-1)^T,\ \ (-1,k-1,k-1)^T.\]
We omit the details, noting that
it is more convenient for us to work subsequently with the inexact estimates
of (\ref{eq4.9}) than with the exact expressions 
thus obtained.$\hfill\square$ \end{rem} 

\begin{cor}\label{cor4.4}
There exists a constant $D_k$, depending only on $k$, such that
\begin{equation}\label{eq4.18}
|\nu_n(\sigma_1,\tau_1) - \nu_n(\sigma_2,\tau_2)|\le D_k
\end{equation}
whenever $\sigma_1,\sigma_2,\tau_1,\tau_2$ are non--empty subsets of $S_k$
satisfying
\begin{equation}\label{eq4.19}
|\sigma_1| = |\sigma_2|,\quad |\tau_1|  = |\tau_2|.
\end{equation}
\end{cor}

\begin{proof}
First suppose that these are one point sets. Then each $\nu_n(\cdot,\cdot)$ is
either $\alpha_n,\beta_n$ or $\gamma_n$, so the difference is estimated by
$2C_k$, using (\ref{eq4.9}). For the general case, (\ref{eq4.1}) shows that
this difference can be realized as a sum of $|\sigma_1| \ |\tau_1|$
differences for one point sets. This gives the estimate
\begin{equation}\label{eq4.20}
|\nu_n(\sigma_1,\tau_1) - \nu_n(\sigma_2,\tau_2)| \le 2|\sigma_1| \ |\tau_1|
C_k.
\end{equation}
Let $D_k$ be the largest possible right hand side in (\ref{eq4.20}), which is
$8k^2C_k$.
\end{proof}

\begin{lem}\label{lem4.5}
Let $x = x_\ell\ldots x_1$ and $y = y_1\ldots y_m$ be words in $\bb F_k$, and
let $\mu(r,s,n;x,y)$ be the number of reduced words in the product $xw_ny$
which result from $r$ cancellations on the left and $s$ cancellations on the
right. Then there exist subsets $\sigma_r(x)$, $\tau_s(y)$ of $S_k$, whose
cardinalities depend only on $r$ and $s$ respectively, such that
\begin{equation}\label{eq4.21}
\mu(r,s,n;x,y) = \nu_{n-r-s}(\sigma_r(x), \tau_s(y)),
\end{equation}
for $n\ge \ell + m+2$, ~ $0 \le r\le\ell$, ~ $0 \le s\le m$.
\end{lem}

\begin{proof}
We first define $\sigma_r(x)$ and $\tau_s(y)$, for $0 \le r\le \ell$, $0\le s
\le m$, and then verify that they have the required properties. Let
\begin{equation}\label{eq4.22}
\sigma_0(x) = S_k\backslash\{x^{-1}_1\}, ~~\sigma_\ell(x) =
S_k\backslash\{x_\ell\},~~ \tau_0(y) = S_k\backslash\{y^{-1}_1\}, ~~\tau_m(y) =
S_k\backslash\{y_m\}.
\end{equation}
For $0 < r < \ell$ and $0 < s < m$, let
\begin{equation}\label{eq4.23}
\sigma_r(x) = S_k\backslash\{x^{-1}_{r+1}, x_r\}, ~~\tau_s(y) = S_k
\backslash\{y^{-1}_{s+1}, y_s\}.
\end{equation}
Note that the pairs $\{x^{-1}_{r+1}, x_r\}$ and $\{y^{-1}_{s+1},y_s\}$ are
distinct, since otherwise cancellations would occur in $x$ or $y$.  Thus the
cardinalities of the sets in (\ref{eq4.22}) and (\ref{eq4.23}) are
respectively $2k-1$ and $2k-2$, and they depend only on $r$ and $s$. If the
reduced word
\begin{equation}\label{eq4.24}
x_\ell \ldots x_{r+1}vy_{s+1}\ldots y_m
\end{equation}
results from one of these cancellations, then the first letter of $v$ must not
be $x^{-1}_{r+1}$ (if $r=\ell$, then this constraint disappears). The original
reduced word in $w_n$ which canceled to this was
\begin{equation}\label{eq4.25}
x^{-1}_1 \ldots x^{-1}_r v y^{-1}_s \ldots y^{-1}_1
\end{equation}
which requires the first letter of $v$ to be different from $x_r$ (this
constraint disappears for $r=0$). Thus the first letter of $v$ must lie in
$\sigma_r(x)$. Conversely, such a $v$ allows exactly $r$ cancellations on the
left in (\ref{eq4.25}). A similar analysis on the right shows that a word of
the form (\ref{eq4.25}) is both reduced and allows the correct number of
cancellations precisely when $v$ starts with an element of $\sigma_r(x)$ and
ends with an element of $\tau_s(y)$. This proves the result.
\end{proof}\newpage

\section{Asymptotic homomorphisms}\label{sec5}

\setcounter{equation}{0}

\indent

We now apply the results of the previous section to show that $\bb E_{\cl B}$
is an asymptotic homomorphism.

\begin{thm}\label{thm5.1}
Let $k\ge 2$ and let $\cl B$ be the Laplacian masa in $VN(\bb F_k)$. Then $\bb
E_{\cl B}$ is an asymptotic homomorphism.
\end{thm}

\begin{proof}
Let $x$ and $y$ in $\bb F_k$ be fixed words of lengths $\ell$ and $m$
respectively. We require that $\ell,m\ge 1$ since otherwise the inequality
(\ref{eq5.10}), for which we are aiming, is trivial. Let $z$ be an arbitrary
word of length $\ell$, and suppose throughout that $n\ge \ell + m + 2$.

If a word $v$ has length $p$ then it is orthogonal to $w_n$ for $n\ne p$.
From this  and (\ref{eq2.7}), it follows that
\begin{equation}\label{eq5.1}
\bb E_{\cl B}(v) = w_p\|w_p\|^{-2}_2.
\end{equation}
Thus, with the notation of Lemma \ref{lem4.5},
\begin{equation}\label{eq5.2}
\bb E_{\cl B}(xw_ny) = \sum^\ell_{r=0} \sum^m_{s=0} \nu_{n-r-s}
(\sigma_r(x), \tau_s(y)) w_{n+\ell+m-2(r+s)}\|w_{n+\ell+m-2(r+s)}\|^{-2}_2.
\end{equation}
By Corollary \ref{cor4.4} and Lemma \ref{lem4.5}
\begin{equation}\label{eq5.3}
|\nu_p(\sigma_r(x), \tau_s(y)) - \nu_p(\sigma_r(z), \tau_s(y)|\le D_k
\end{equation}
for $p\ge 2$. Thus there exist constants $\lambda_{r,s}$, uniformly bounded by
$D_k$, such that 
\begin{equation}\label{eq5.4}
\bb E_{\cl B}(xw_ny) - \bb E_{\cl B}(zw_ny) = \sum^\ell_{r=0} \sum^m_{s=0}
\lambda_{r,s}w_{n+\ell+m-2(r+s)}\|w_{n+\ell+m-2(r+s)}\|^{-2}_2.
\end{equation}
Now $\|w_p\|^2_2$ is the number of terms in $w_p$, so
\begin{equation}\label{eq5.5}
\|w_p\|^2_2 = 2k(2k-1)^{p-1} = (2k-1)^{p-n} \|w_n\|^2_2.
\end{equation}
Thus
\begin{equation}\label{eq5.6}
\|w_{n+\ell+m-2(r+s)}\|^2_2 = (2k-1)^{\ell+m-2(r+s)}\|w_n\|^2_2,
\end{equation}
so
\begin{equation}\label{eq5.7}
\|w_{n+\ell+m-2(r+s)}\|^2_2 \ge (2k-1)^{-(\ell+m)} \|w_n\|^2_2,
\end{equation}
for $0 \le r\le \ell, \ 0\le s\le m$. This last inequality and (\ref{eq5.4})
then give the estimate
\begin{equation}\label{eq5.8}
\|\bb E_{\cl B}(xw_ny) - \bb E_{\cl B}(zw_ny)\|_2 \le (\ell +1)
(m+1)D_k(2k-1)^{(\ell+m)/2} \|w_n\|^{-1}_2.
\end{equation}
Let $H_{\ell,m,k}$ be the constant on the right hand side of (\ref{eq5.8}). If
we sum (\ref{eq5.8}) over all words $z$ of length $\ell$ (of which there are
$\|w_\ell\|^2_2$), we obtain
\begin{equation}\label{eq5.9}
\|~\|w_\ell\|^2_2 \,\bb E_{\cl B}(xw_ny) - \bb E_{\cl B}(w_\ell w_ny)\|_2 \le
\|w_\ell\|^2_2 H_{\ell,m,k} \|w_n\|^{-1}_2.
\end{equation}
Since $\bb E_{\cl B}(x) = w_\ell\|w_\ell\|^{-2}_2$ and $\bb E_{\cl B}(w_\ell
w_ny) = w_\ell w_n\bb E_{\cl B}(y)$, (\ref{eq5.9}) implies that
\begin{equation}\label{eq5.10}
\|\bb E_{\cl B}(xw_ny) - \bb E_{\cl B}(x) \bb E_{\cl B}(y)w_n\|_2 \le
H_{\ell,m,k} \|w_n\|^{-1}_2.
\end{equation}

If we let $v_n = w_n \|w_n\|^{-1}_2$, then the terms of the series in
(\ref{eq3.1}) are bounded by $H^2_{\ell,n,k}\|w_n\|^{-4}_2$, and the series is
clearly summable. The result now follows from Theorem~\ref{thm3.1}.
\end{proof}

\begin{thm}\label{thm5.2}
Fix $k\ge 2$. Then the Laplacian masa $\cl  B$ in $\bb F_k$ is strongly
singular, and has Popa invariant $\delta(\cl B)=1$.
\end{thm}

\begin{proof}
Apply Theorem \ref{thm5.1} and \cite[Theorem~6.7]{SS2001}. 
\end{proof}

\begin{rem}\label{rem5.3}
Although we have referred throughout to $\cl B$ as a masa, we have never made
use of this fact. Thus our results give new proofs that $\cl B$ is a masa and
is singular, \cite{Py,Ra}.$\hfill\square$
\end{rem}

We now extend our results to the free product of a finite number of countable 
discrete groups $G_i$, $1\leq i \leq m$, $m \geq 2$. We let $G$ denote
$G_1*G_2*\ldots *G_m$, and we fix elements $g_i \in G_i$, $1 \leq i \leq k$.
We assume that each $g_i$ has infinite order, which is possible in many, but
not all, such groups. We denote by $\cl B$ the abelian subalgebra of $VN(G)$
generated by the self-adjoint element
\[h=\sum_{i=1}^k (g_i+g_i^{-1}).\]
 \begin{thm}\label{thm5.4}
With the above notation, the conditional expectation $\bb E_{\cl B}$ of 
$VN(G)$ onto $\cl B$ is an asymptotic homomorphism, $\cl B$ is a strongly
singular masa, and $\delta(\cl B)=1$.
\end{thm}

\begin{proof}
The elements $g_1,\ldots ,g_k$ generate a copy of ${\bb F}_k$ inside $G$, and
we write $\cl N$ for the resulting subfactor of $VN(G)$. As before, we let
$w_n$ be the sum of all words of length $n$ in ${\bb F}_k$, and we note that 
$\{v_n=w_n/\|w_n\|_2\}_{n=0}^{\infty}$ is an orthonormal basis for $L^2(\cl
B)$. The result will follow from Theorem \ref{thm3.1} if we can show that 
\begin{equation}\label{eq5.11}
\sum^\infty_{n=0} \|\bb E_{\cl B}(xv_ny) - \bb E_{\cl B}(x) \bb E_{\cl B}(y)
v_n\|^2_2 < \infty
\end{equation}
for all $x,y \in G$.

If $x \in {\bb F}_k$, $y \in G$, then 
\begin{equation}\label{eq5.12}
\bb E_{\cl B}(xv_ny)=\bb E_{\cl B}(\bb E_{\cl N}(xv_ny))=
\bb E_{\cl B}(xv_n\bb E_{\cl N}(y)),
\end{equation}
using modularity of $\bb E_{\cl N}$. If $y \in G\backslash {\bb F}_k$, then
$\bb E_{\cl N}(y)=0$, so all terms in (\ref{eq5.11}) vanish. If $y \in {\bb
F}_k$, then the argument of Theorem \ref{thm5.1} shows the validity of
(\ref{eq5.11}). A similar analysis holds if we begin by supposing that $y \in
{\bb F}_k$. We may thus assume that $x,y \in G\backslash {\bb F}_k$. In this
case $\bb E_{\cl B}(x)=0$, so (\ref{eq5.11}) becomes
\begin{equation}\label{eq5.13}
\sum^\infty_{n=0} \|\bb E_{\cl B}(xv_ny)\|^2_2 < \infty .
\end{equation}
The elements $x$ and $y$ may be written as products
\begin{equation}\label{eq5.14}
x=x_{\ell}\ldots x_1,\ \ \ y=y_1\ldots y_q,
\end{equation}
where $x_i \in G_{s_i}$, $y_j \in G_{t_j}$, and $s_i \neq s_{i+1}$, $t_j \neq
t_{j+1}$. Let
\[\chi_n=\{u \in {\bb F}_k:\,|u|=n,\ xuy \in {\bb F}_k\},\]
where $|u|$ is the length of $u$ in terms of $g_1,\ldots ,g_k \in {\bb
F}_k$. We will show that the cardinality $|\chi_n|$ is at most $(2n+1)(n+1)$.
Consider a word $u \in {\bb F}_k$ of length $n$ and let 
\[u=u_1\ldots u_r\] be written in the free product $G$ with each $u_j$ a
non--zero power of some $g_{i_j}$, and $i_j \neq i_{j+1}$ for all $j$. Note
that $r \leq n$. We observe that the product $xuy$ is in ${\bb F}_k$ if and
only if one of the following conditions is satisfied:\newline
(1) there exists $p$, $0 \leq p \leq r$, such that $u_i=x_i^{-1}$ for $1 \leq
i \leq p$, and $u_{r-i+1}=y_i^{-1}$ for $1 \leq i \leq r-p$ (with obvious
modifications if $p=0$ or $p=r$), and the element
\begin{equation}\label{eq5.15}
xuy=x_{\ell}\ldots x_{p+1}y_{r-p+1}\ldots y_q
\end{equation}
is in ${\bb F}_k$;\newline
(2) there exists $p$, $0 \leq p \leq r$, such that $u_i=x_i^{-1}$ for $1 \leq
i \leq p-1$, and $u_{r-i+1}=y_i^{-1}$ for $1 \leq i \leq r-p$ (with obvious
modifications if $p=0$ or $p=r$), and the element
\begin{equation}\label{eq5.16}
xuy=x_{\ell}\ldots x_{p}u_py_{r-p+1}\ldots y_q
\end{equation}
is in ${\bb F}_k$.

In (\ref{eq5.15}), $x_{p+1}$ and $y_{r-p+1}$ are in the same $G_i$ since
further cancellation must occur, while in (\ref{eq5.16}) the same conclusion
applies to the elements $x_p$, $u_p$, and $y_{r-p+1}$. In (1), there are at
most $n+1$ values of $p$, and thus at most the same number of possibilities
for $u$. In (2), $p$ takes up to $n+1$ values, and the remaining term $u_p$ 
is a non--zero power $g_i^s$ for some $i$ and some $s$ with $1\leq |s| \leq n$.
This gives an upper bound of $2n(n+1)$ possibilities. The two cases combine to
give an upper estimate (which could undoubtedly be lowered by a more detailed
analysis) of $(n+1)(2n+1)$ for $|\chi_n|$.

It now follows that 
\begin{equation}\label{eq5.17}
\|\bb E_{\cl B}(xw_ny)\|_2 \leq |\chi_n| \leq (n+1)(2n+1),
\end{equation}
since $\|{\bb E}_{\cl B}(xuy)\|_2 \leq 1$ for each $u \in \chi_n$. Thus the
terms in (\ref{eq5.13}) are dominated by $(n+1)^2(2n+1)^2\|w_n\|_2^{-2}$, so
convergence is guaranteed by (\ref{eq2.6}). This proves that $\bb E_{\cl B}$
is an asymptotic homomorphism, and the other statements now follow from 
\cite[Theorem 6.7]{SS2001}.
\end{proof}

\begin{rem}\label{rem5.5}
The conclusions of Theorem \ref{thm5.4} would remain valid were we to replace
each $g_i$ by some non--zero power $g_i^{t_i}$, and $h$ by
\[\sum_{i=1}^k (g_i^{t_i}+g_i^{-t_i}),\]
since $\{g_i^{t_i}\}_{i=1}^k$ also generates a copy of ${\bb F}_k$
in $G$. $\hfill\square$
\end{rem}

\end{document}